\documentclass[12pt]{article}
\usepackage{amssymb}
\usepackage{graphics}
\usepackage{latexsym}
\usepackage{amsfonts}
\usepackage{amscd}
\usepackage{makeidx}
\usepackage{epsfig}
\usepackage{amsmath}
\usepackage{color}

\newtheorem{theorem}{Theorem}[section]

\newtheorem{proposition}[theorem]{Proposition}

\input epsf

\begin{document}

\title{On a Lie algebraic structure associated with a non-linear
dynamical system }
\author{J.R. Guzm\'{a}n \thanks{%
Econom\'{\i}a Aplicada. Instituto de Investigaciones Econ\'{o}micas.
Universidad Nacional Aut\'{o}noma de M\'{e}xico. e-mail: \ jrg@unam.mx. \
This article was carried out with the support of a grant from from
UNAM-DGAPA during \ the sabbatical \ year in the PUIMECI of the UACh.}}
\maketitle

\begin{abstract}
A family of Lie algebras of minimal dimension associated with vector fields
that define a non-linear dynamical system is calculated. These Lie algebras
contain the Heinsenberg algebra.

An element that distinguishes these vector fields is called
evapotranspiration function. This function can be \ calculated solving
equations in partial derivatives that arise in determining the Heinsenberg
algebra.

Using Kozsul homology for this  Lie algebras,  Euler characteristic  is
calculated.
\end{abstract}

\section{Introduction}

The system that appears in \cite{Blasco} provided by the differential equation
system below, is considered:

\begin{eqnarray*}\label{sistema}
\left( 
\begin{array}{c}
\stackrel{\cdot }{x} \\ 
\stackrel{\cdot }{y} \\ 
\stackrel{\cdot }{z}%
\end{array}%
\right) &=&\left( 
\begin{array}{c}
\alpha _{1} \\ 
\alpha _{2} \\ 
\alpha _{3}%
\end{array}%
\right) +\left( 
\begin{array}{c}
\beta _{11}x+\beta _{12}E\left( x,y\right) \\ 
\beta _{22}y+\beta _{22}^{\prime }E\left( x,y\right) +\beta _{13}z \\ 
\beta _{32}y+\beta _{33}z%
\end{array}%
\right) \\
&&+\alpha \left( 
\begin{array}{c}
\gamma _{11} \\ 
\gamma _{21} \\ 
0%
\end{array}%
\right) +f\left( 
\begin{array}{c}
\gamma _{12} \\ 
\gamma _{22} \\ 
0%
\end{array}%
\right) +w\left( 
\begin{array}{c}
\gamma _{13} \\ 
\gamma _{23} \\ 
0%
\end{array}%
\right) .
\end{eqnarray*}

This dynamical system is a simplification of the original system; $x,y$ and $z$
denote $x_{i},T_{i},T_{m};$ are \ the state variables that represent
relative humidity, air temperature and thermal mass temperature. It is a
dynamic system that has been formulated using thermodynamic principles.

In this system,\eqref{sistema} all the thermodynamic functions the authors propose have
been substituted. The letters $\alpha ,\beta ,\gamma ,$and their sub-indices
depend on original parameters.

$\alpha ,$ $f,$ $w$ are control variables that represent the window opening
angle, nebulization system intensity and heating system intensity, respectively.
Although all of this work is based on what is itself a dynamic control system, these
variables are not used in the subsequent.

$\alpha ,$ $\beta ,$ $\gamma $ and sub-indices are dynamic system
parameters; the signs of the parameters are $\alpha _{i}>0,$ for $i=1,2,3;$ $%
\beta _{11},\beta _{22},\beta _{33},\beta _{22}^{\prime }<0,$ $\beta
_{12},\beta _{13},\beta _{32}>0,$ $\gamma _{11},\gamma _{12},\gamma _{23}>0,$
$\gamma _{21},\gamma _{22}<0,$ $\ \gamma _{13}=0.$

Although in the original article \cite{Blasco} appears $\gamma _{13}$ as null, is
considered here $\ \gamma _{13}>0$; for consistency with the units specified
in the dynamic system should be $\gamma _{13}[kg_{air}^{-1}]$ units. This
change represents a \ heating on the humidity variable.

It is worth pointing out that in the original dynamic system,
$E\left( x,y\right) $ is the evapotranspiration function,
only appears as a function of \ $y$. In this article, evapotranspiration $E$
is generalized and made to depend on $x.$

Crops in a
greenhouse are subject to an internal dynamic that creates
evapotranspiration $E[kgH_{2}Os^{-1}]$;
this function provides a number of
estimation proposals, depending on each author. Evapotranspiration research
is an active research field, as testified to by \ [4], for example. Function 
$E$, is functioning as an integral part of the model, and dynamic richness
is lent to the system in question. This definition can be stipulated since
evapotranspiration estimation is an open problem.

In this article, a family of functions of evapotranspiration depending of
the $\gamma _{ij}$ \ parameters, which naturally arises from the algebraic
structure, can be calculated to define Lie algebra structure on the dynamic
system into consideration.

\noindent
More specifically. This dynamic system is of the form
$$\stackrel{\cdot }{X}%
=f_{0}(X)+\underset{i=1}{\sum}u_{i}f_{i}(X),
$$

\noindent
often investigated in control theory; $X$ is a curve in $\mathbb{R}^{n}$, \ $f_{i}$%
\'{}%
s are sufficiently differentiable functions, $u_{i}$%
\'{}%
s are the control variables. In the case that concerns,
$$X=\left( 
\begin{array}{c}
x \\ 
y \\ 
z%
\end{array}%
\right) , f_{0}(X)=\left( 
\begin{array}{c}
\alpha _{1} \\ 
\alpha _{2} \\ 
\alpha _{3}%
\end{array}%
\right) +\left( 
\begin{array}{c}
\beta _{11}x+\beta _{12}E\left( x,y\right) \\ 
\beta _{22}y+\beta _{22}^{\prime }E\left( x,y\right) +\beta _{13}z \\ 
\beta _{32}y+\beta _{33}z%
\end{array}%
\right) $$ $$ and \ \ f_{i}(X)=\left( 
\begin{array}{c}
\gamma _{1i} \\ 
\gamma _{2i} \\ 
0%
\end{array}%
\right) ,i=1,2,3.
$$ 
In addition $u_{1}=\alpha ,$ $u_{2}=f,$ $u_{3}=w$. With
these functions its posible define vector fields that are used to elucidate
this algebraic structure of this dynamic system, (see below, in \ section \
\ of proof of proposition 1). The Lie bracket is the usual for vector fields.

From these calculations \ and contained in the infinite dimensional family
of \ Lie algebras it is had a minimal \ dimension family of Lie algebras in
which is contained the \ Heinsenberg algebra.

As it is well known, this algebra plays an important role in investigations
of energy of quantum systems. This algebra began to be considered by the
physicist Werner Heisenberg in relation to quantic oscillators.
For current research on this topic see \cite{Morton}. For a pleasant way
of explanation about this, see \cite{Baez-Morton} . Another example is in its application to
gravitational physics of \cite{Tkachuk}. The Heinsenberg algebras play a central role
in the research of new ways of understanding the gravitational field.

By completeness is calculated the homology \ of the family of Lie algebras \ 
$\mathfrak{a}_{E}.$The homology is coincident with the results known to the
case of the $\mathfrak{h}_{3}$ algebra. The homology of these algebras to
been objet of a deep study. The pure Lie algebra homology has been
considered in work starting with \cite{Santaroubane} up to \cite{Cairns}, for example.

In general, there are investigations that are a central part of the homology
and dynamical systems. For a review of the most important theorems, see \cite{Sanchez}.

The relevance of these results is that when conditions searching for the
definition of the dynamic system underlying Lie algebra, a function of
evapotranspiration is the solution of partial differential equations that
arise in applying algebraic Lie tools.

\section{Main results}

The following are the main results.

\begin{proposition}
The \ family of Lie algebras $\mathfrak{a}_{E},$ \
of minimal dimension \ 7, is specified by conditioning that $\Delta
_{ij}E\left( x,y\right) =c_{ij}$; and $\Delta _{ijk}E\left( x,y\right) =0$
and $\Lambda _{ijk}E\left( x,y\right) =0.$

$\mathfrak{a}_{E}$ is the Heinsenberg algebra if $\Delta _{ij}E\left(
x,y\right) $ does not equal 0 for $i=j$ and \ $\Delta _{ij}E\left(
x,y\right) =0$, for $i\neq j$; \ where

\begin{eqnarray*}
\Delta _{ij}E\left( x,y\right) &=&\gamma _{1i}\gamma _{1j}\frac{\partial
^{2}E}{\partial x^{2}}+2\left( \gamma _{1i}\gamma _{2j}+\gamma _{1j}\gamma
_{2i}\right) \frac{\partial ^{2}E}{\partial x\partial y} \\
&&+\gamma _{2i}\gamma _{2j}\frac{\partial ^{2}E}{\partial y^{2}};
\end{eqnarray*}

\noindent
and

\begin{eqnarray*}
\Delta _{ijk}E\left( x,y\right) &=&\gamma _{1i}\gamma _{1j}\gamma _{1k}\frac{%
\partial ^{3}E}{\partial x^{3}}+\gamma _{2i}\gamma _{2j}\gamma _{2k}\frac{%
\partial ^{3}E}{\partial y^{3}}+ \\
&&\frac{\partial ^{3}E}{\partial x^{2}\partial y}\underset{\left( \sigma
_{1},\sigma _{2},\sigma _{3}\right) \in \mathfrak{m}_{3}\left( 4\right) }{%
\sum }\gamma _{\sigma _{1}i}\gamma _{\sigma _{2}j}\gamma _{\sigma _{3}k}+ \\
&&\frac{\partial ^{3}E}{\partial x\partial y^{2}}\underset{\left( \sigma
_{1},\sigma _{2},\sigma _{3}\right) \in \mathfrak{m}_{3}\left( 5\right) }{%
\sum }\gamma _{\sigma _{1}i}\gamma _{\sigma _{2}j}\gamma _{\sigma _{3}k};
\end{eqnarray*}

\noindent
where $\mathfrak{m}_{k}\left( n\right) $ is the magic square of $k\times k$
with magic sum $n$. The notation $\sigma \in \mathfrak{m}_{k}\left( n\right) 
$ means taking line $\sigma $ from $\mathfrak{m}_{k}\left( n\right) $.

\noindent
And

$\Lambda _{ijk}E\left( x,y\right) $ $=$

$\beta _{12}\gamma _{1i}\gamma _{1j}\gamma _{1k}\left( \frac{\partial ^{2}E}{%
\partial x^{2}}\right) ^{2}+\beta _{22}^{\prime }\gamma _{2i}\gamma
_{2j}\gamma _{2k}\left( \frac{\partial ^{2}E}{\partial y^{2}}\right) ^{2}$

$+$($\beta _{12}\gamma _{2i}\gamma _{2j}\gamma _{1k}+\beta _{22}^{\prime
}\gamma _{1i}\gamma _{1j}\gamma _{2k}$)$\frac{\partial ^{2}E}{\partial x^{2}}%
\frac{\partial ^{2}E}{\partial y^{2}}$

$+$($\beta _{12}\gamma _{1i}\gamma _{2j}\gamma _{2k}+\beta _{12}\gamma
_{2i}\gamma _{1j}\gamma _{2k}+\beta _{22}^{\prime }\gamma _{1i}\gamma
_{2j}\gamma _{1k}+\beta _{22}^{\prime }\gamma _{2i}\gamma _{1j}\gamma _{1k}$)%
$\left( \frac{\partial ^{2}E}{\partial x\partial y}\right) ^{2}$

$+$($(\beta _{12}\gamma _{1i}\gamma _{1j}\gamma _{2k}+\beta _{12}\gamma
_{1i}\gamma _{2j}\gamma _{1k}+\beta _{12}\gamma _{2i}\gamma _{1j}\gamma
_{1k}+\beta _{22}^{\prime }\gamma _{1i}\gamma _{1j}\gamma _{1k})\frac{%
\partial ^{2}E}{\partial x^{2}}$

$+$ \ $(\beta _{12}\gamma _{2i}\gamma _{2j}\gamma _{2k}+\beta _{22}^{\prime
}\gamma _{1i}\gamma _{2j}\gamma _{2k}+\beta _{22}^{\prime }\gamma
_{2i}\gamma _{1j}\gamma _{2k}+\beta _{22}^{\prime }\gamma _{2i}\gamma
_{2j}\gamma _{1k})\frac{\partial ^{2}E}{\partial y^{2}}$)$\frac{\partial
^{2}E}{\partial x\partial y}$

$-$($+\beta _{12}\gamma _{1i}\gamma _{1j}\gamma _{1k}\frac{\partial E}{%
\partial x}+\beta _{12}\gamma _{1i}\gamma _{1j}\gamma _{2k}\frac{\partial E}{%
\partial y}+\beta _{11}\gamma _{1i}\gamma _{1j}\gamma _{1k}$)$\frac{\partial
^{3}E}{\partial x^{3}}$

$-$($\beta _{11}\gamma _{1i}\gamma _{2j}\gamma _{1k}+\beta _{11}\gamma
_{2i}\gamma _{1j}\gamma _{1k}+\beta _{22}\gamma _{1i}\gamma _{1j}\gamma
_{2k} $

$-$($\beta _{12}\gamma _{1i}\gamma _{2j}\gamma _{1k}+\beta _{12}\gamma
_{2i}\gamma _{1j}\gamma _{1k}+\beta _{22}^{\prime }\gamma _{1i}\gamma
_{1j}\gamma _{1k}$)$\frac{\partial E}{\partial x}$

$-$($\beta _{12}\gamma _{1i}\gamma _{2j}\gamma _{2k}+\beta _{12}\gamma
_{2i}\gamma _{1j}\gamma _{2k}+\beta _{22}^{\prime }\gamma _{1i}\gamma
_{1j}\gamma _{2k}$)$\frac{\partial E}{\partial y}$)$\frac{\partial ^{3}E}{%
\partial x^{2}\partial y}$

$-$($\beta _{11}\gamma _{2i}\gamma _{2j}\gamma _{1k}+\beta _{22}\gamma
_{1i}\gamma _{2j}\gamma _{2k}+\beta _{22}\gamma _{2i}\gamma _{1j}\gamma
_{2k} $

$-$($\beta _{12}\gamma _{2i}\gamma _{2j}\gamma _{1k}+\beta _{22}^{\prime
}\gamma _{1i}\gamma _{2j}\gamma _{1k}+\beta _{22}^{\prime }\gamma
_{2i}\gamma _{1j}\gamma _{1k}$)$\frac{\partial E}{\partial x}$

$-$($\beta _{12}\gamma _{2i}\gamma _{2j}\gamma _{2k}+\beta _{22}^{\prime
}\gamma _{1i}\gamma _{2j}\gamma _{2k}+\beta _{22}^{\prime }\gamma
_{2i}\gamma _{1j}\gamma _{2k}$)$\frac{\partial E}{\partial y}$)$\frac{%
\partial ^{3}E}{\partial x\partial y^{2}}$

$-$($\beta _{22}^{\prime }\gamma _{2i}\gamma _{2j}\gamma _{1k}\frac{\partial
E}{\partial x}+\beta _{22}^{\prime }\gamma _{2i}\gamma _{2j}\gamma _{2k}%
\frac{\partial E}{\partial y}+\beta _{22}\gamma _{2i}\gamma _{2j}\gamma
_{2k} $)$\frac{\partial ^{3}E}{\partial y^{3}}$
\end{proposition}

\medskip

\begin{proposition}
\begin{equation*}
E_{i}(x,y)=\frac{1}{C_{4}+C_{3}\tanh (\frac{C_{1}\gamma _{1i}-C_{2}\gamma
_{2i}x+C_{2}\gamma _{1i}y}{\gamma _{1i}})},\text{ }
\end{equation*}

\noindent
the $C_{k}^{\prime }s$ \ are arbitrary constants.
\end{proposition}

\begin{proposition} 
The Euler Characteristic of the \ generated
subalgebras $\mathfrak{a}_{E}$ for the dynamic systems is $\mathcal{X}\left( 
\mathfrak{a}_{E}\right) =0.$
\end{proposition}
\section{Proof of Proposition 2.1}

The proof of proposition 1 is as follows.

\noindent
To elucidate the structure of Lie algebras is must associated vector fields
the dynamic system mentioned above.

\noindent
This dynamic system \eqref{sistema} is characterized by vector fields

\begin{eqnarray*}
f_{0} &=&\left( \alpha _{1}+\beta _{11}x+\beta _{12}E\left( x,y\right)
\right) \frac{\partial }{\partial x}+ \\
&&\left( \alpha _{2}+\beta _{22}y+\beta _{22}^{\prime }E\left( x,y\right)
+\beta _{13}z\right) \frac{\partial }{\partial y}+ \\
&&\left( \alpha _{3}+\beta _{32}y+\beta _{33}z\right) \frac{\partial }{%
\partial z}, \\
f_{i} &=&\gamma _{1i}\frac{\partial }{\partial x}+\gamma _{2i}\frac{\partial 
}{\partial y},\text{ for }i=1,2,3.
\end{eqnarray*}

\noindent
The vector field $B$ is defined as

\begin{equation*}
B=-\beta _{12}\frac{\partial }{\partial x}-\beta _{22}^{\prime }\frac{%
\partial }{\partial y},
\end{equation*}

\noindent
called here the air vector; $\beta _{12}$ features $kg_{air}^{-1},$ physical
units, and is a dimensionless constant of the contained air volume. $\beta
_{22}^{\prime }$ features $joule\cdot kg_{air}^{-1},$ units, i.e., energy
contained in the air.

$\mathfrak{a}_{E}$ sub-algebras can be defined using the multiplication
table given by the following relationships

\begin{eqnarray*}
\left[ f_{i},\left[ f_{0},f_{j}\right] \right] &=&\Delta _{ij}E\left(
x,y\right) B,\text{ where }i,j=1,2,3 \\
\left[ f_{i},f_{j}\right] &=&0 \\
\left[ \left[ f_{0},f_{i}\right] ,\left[ f_{0},f_{j}\right] \right] &=&0,
\end{eqnarray*}

\noindent
where

\begin{eqnarray*}
\left[ f_{0},f_{j}\right] &=&-\left( \gamma _{1j}\beta _{11}+\gamma
_{1j}\beta _{12}\frac{\partial E}{\partial x}+\gamma _{2j}\beta _{12}\frac{%
\partial E}{\partial y}\right) \frac{\partial }{\partial x} \\
&&-\left( \gamma _{2j}\beta _{22}+\gamma _{1j}\beta _{22}^{\prime }\frac{%
\partial E}{\partial x}+\gamma _{2j}\beta _{22}^{\prime }\frac{\partial E}{%
\partial y}\right) \frac{\partial }{\partial y} \\
&&-\gamma _{2j}\beta _{32}\frac{\partial }{\partial z}.
\end{eqnarray*}

The relationship $\Delta _{ij}E\left( x,y\right) =\Delta _{ji}E\left(
x,y\right) $ is satisfied.

\noindent
In order that the non-trivial relationship $\left[ f_{i},\left[ f_{0},f_{j}%
\right] \right] =\Delta _{ij}E\left( x,y\right) B$ define a multiplication
table, the first condition which must be met is that $\left[ f_{k},\Delta
_{ij}E\left( x,y\right) B\right] =\Delta _{ijk}E\left( x,y\right) B=0$. And
the second condition that must be met is that

\begin{equation*}
\left[ \left[ f_{0},f_{k}\right] ,\Delta _{ij}E\left( x,y\right) B\right]
=\Lambda _{ijk}E\left( x,y\right) B=0.
\end{equation*}

In sum, in order for Lie algebra $\mathfrak{a}_{E}$ be defined it must have $%
\Delta _{ij}E\left( x,y\right) =c_{ij};$ where $c_{ij}$ are constants and as
well, $\Delta _{ijk}E\left( x,y\right) =0$ \ \ and $\ \Lambda _{ijk}E\left(
x,y\right) =0.$

It should be noted that the $B$ continues to appear when other Lie
parentheses are calculated, like $\left[ f_{k1},\left[ f_{k2},...,\left[
f_{kl},\Delta _{ij}E\left( x,y\right) B\right] \right] \right] =(\cdot )B$,
where $(\cdot )$ is an expression that contains partial derivatives of $E$
of some determined order. And the sums contained in these expressions are
defined in terms of superior order magic squares and greater magic sums.

The fact that the vector field $B$ appears frequently in the calculations as
a factor in the final \ result, was decisive for \ searching a Lie algebra
structure

\section{Proof of Proposition 2.2}

To solve the system of partial differential equations

\begin{eqnarray*}
\Delta _{ij}E\left( x,y\right) &=&\left\{ 
\begin{array}{c}
1\text{ }for\text{ }i=j \\ 
0\text{ }for\text{ }i\neq j%
\end{array}%
\right. \\
\Delta _{ijk}E\left( x,y\right) &=&0 \\
\Lambda _{ijk}E\left( x,y\right) &=&0,
\end{eqnarray*}

\noindent
can have the function of evapotraspiracion $E$.

In general the conditions that define $\mathfrak{a}_{E}$ sub-algebras and $%
\mathfrak{h}_{3}$ sub-algebra are non linear differential equation systems
in partial derivatives with dependent variable $E$. Nevertheless within the
complexity that searching for $E$ solutions in these systems can represent,
it is possible to find general evapotranspiration functions with adequate
boundary conditions.

For the system of partial differential equations that define $\mathfrak{h}%
_{3}$, the system of equations can be resolved separately, in three
subsystems.

The first class of equations \ $\Delta _{ij}E\left( x,y\right) =\left\{ 
\begin{array}{c}
1\text{ }for\text{ }i=j \\ 
0\text{ }for\text{ }i\neq j%
\end{array}%
\right.$ has the general solution.

For $i\neq j$

\begin{equation*}
E_{ij}(x,y)=F_{1}(-\gamma _{2j}x+\gamma _{1j}y)+F_{2}(-\gamma _{2i}x+\gamma
_{1i}y).\text{ }
\end{equation*}

\noindent
For $i=j$

\begin{equation*}
E_{ij}(x,y)=F_{1}(-\gamma _{2j}x+\gamma _{1j}y)+F_{2}(-\gamma _{2i}x+\gamma
_{1i}y)+\frac{x^{2}}{2\gamma _{1i}\gamma _{1j}},
\end{equation*}

\noindent
with $F_{1}$, $F$ $_{2}$ arbitrary, sufficiently differentiable functions.

The second class of functions $\Delta _{ijk}E\left( x,y\right) =0$, has the
general solution.

\begin{equation*}
E_{ijk}(x,y)=F_{1}(-\gamma _{2i}x+\gamma _{1i}y)+F_{2}(-\gamma _{2j}x+\gamma
_{1j}y)+F_{3}(-\gamma _{2k}x+\gamma _{1k}y),\text{ }
\end{equation*}

\noindent
with $F_{1},$ $F_{2},$ $F_{3}$ \ arbitrary, sufficiently differentiable
functions.

Because of the voluminous calculations, to resolve the third system of
equations $\Lambda _{ijk}E\left( x,y\right) =0$, for all $i,j,k$ is
possible to use the algebraic package ??? program. The result
of the search for the solution of the system is simple; it is the formula
given in the proposition 2, above enunciated.

Thus, finding functions $F_{i}$'s in the first and second class of solutions
it is possible to calculate the function $E(x,y)$.

\section{Proof of Proposition 2.3}

An important question is the general nature of the dynamic system that the
greenhouse models, along with the conditions that define the Lie
sub-algebras. Is well known that the best invariant to see if the same
sub-algebras are exhibited---given another dynamic system with three state
variables and three control variables by example---is to calculate the
homology of the sub-algebras that are found.

\noindent
The results of the calculation of the homology is as follows.

The $\mathfrak{a}_{E}$ algebra is generated as a \ 7 dimensional vector
space by

\begin{equation*}
\left\{ X_{i}\equiv f_{i},\text{ }Y_{i}\equiv \left[ f_{0},f_{i}\right] ,%
\text{ }Z\equiv B\right\} _{i=1,2,3}.
\end{equation*}

\noindent
Koszul's homology is used on field-$R$, of real numbers, and the trivial
representation \ with the usual defined, $p$-chain to $p-1$-chains border
operator in the alternate $\bigwedge \mathfrak{a}_{E},$ algebra.

\begin{equation*}
\partial _{p}:\bigwedge \mathfrak{a}_{E}\rightarrow \bigwedge \mathfrak{a}%
_{E},
\end{equation*}

\begin{equation*}
\partial _{p}\left( X_{1}\wedge X_{2}\wedge ...\wedge X_{P}\right) =\underset%
{1\leq i<j\leq P}{\sum }\left( -1\right) ^{i+j}\left[ X_{i},X_{j}\right]
\wedge X_{1}\wedge ...\wedge \overset{-}{X}_{i}\wedge ...\wedge \overset{-}{X%
}_{j}\wedge ...\wedge X_{p}.
\end{equation*}

\noindent
The complex to calculate the homology is

\begin{eqnarray*}
&&0\overset{\partial _{8}}{\rightarrow }\underset{\downarrow }{\wedge ^{7}%
\mathfrak{a}_{E}}\overset{\partial _{7}}{\rightarrow }\text{ }\underset{%
\downarrow }{\wedge ^{6}\mathfrak{a}_{E}}\overset{\partial _{6}}{\rightarrow 
}\underset{\downarrow }{\wedge ^{5}\mathfrak{a}_{E}}\overset{\partial _{5}}{%
\rightarrow }\underset{\downarrow }{\wedge ^{4}\mathfrak{a}_{E}}\overset{%
\partial _{4}}{\rightarrow }\underset{\downarrow }{\wedge ^{3}\mathfrak{a}%
_{E}}\overset{\partial _{3}}{\rightarrow }\underset{\downarrow }{\wedge ^{2}%
\mathfrak{a}_{E}}\overset{\partial _{2}}{\rightarrow }\underset{\downarrow }{%
\mathfrak{a}_{E}}\overset{\partial _{1}}{\rightarrow }\text{ }\underset{%
\downarrow }{R}\overset{\partial _{0}}{\rightarrow }0\text{\ \ } \\
\text{\ }0 &\rightarrow &\text{ \ \ \ \ \ \ \ \ \ }R^{7}\ \ \rightarrow 
\text{ \ \ }R^{21}\rightarrow \text{ }R^{35}\rightarrow \text{ }%
R^{35}\rightarrow \text{ \ }R^{21}\rightarrow \text{ \ \ }R^{7}\rightarrow 
\text{ }R\rightarrow \text{ }R\rightarrow 0.
\end{eqnarray*}

\noindent
The matrices associated with $\partial _{k}$ linear transformations that we
will denote as $\partial _{k}^{ij}$ are 
\begin{equation*}
\partial _{7}^{ij}=0\text{, }i=1,...,7;\text{ }j=1.
\end{equation*}

\begin{eqnarray*}
\partial _{6}^{ij} &=&-c_{33},\text{ }i=8;\text{ }j=1 \\
&=&c_{23},\text{ }i=9,12;\text{ }j=1 \\
&=&-c_{13},\text{ }i=10,17;\text{ }j=1 \\
&=&-c_{22},\text{ }i=13;\text{ }j=1 \\
&=&c_{12},\text{ }i=14,18;\text{ }j=1 \\
&=&-c_{11},\text{ }i=19;\text{ }j=1 \\
&=&0\text{ in the other cases of }i,j.
\end{eqnarray*}

Matrices are defined thus

\begin{equation*}
\mathfrak{a}_{2}=\left( 
\begin{array}{ccc}
0 & 0 & 0 \\ 
1 & 0 & 0 \\ 
0 & 1 & 0%
\end{array}%
\right) ,\mathfrak{a}_{0}=\left( 
\begin{array}{ccc}
-1 & 0 & 0 \\ 
0 & 0 & 0 \\ 
0 & 0 & 1%
\end{array}%
\right) ,\mathfrak{a}_{-2}=\left( 
\begin{array}{ccc}
0 & -1 & 0 \\ 
0 & 0 & -1 \\ 
0 & 0 & 0%
\end{array}%
\right).
\end{equation*}

Note that the elements $\mathfrak{a}_{2},\mathfrak{a}_{0},\mathfrak{a}_{-2}$%
, form non-trivial commutators relations [$\mathfrak{a}_{i},\mathfrak{a}_{j}$%
]=$-\mathfrak{a}_{i+j},$ forming a base of some Lie algebra.

Some $\partial _{k}^{ij}$ that follow are expressed by means of $\mathfrak{a}%
_{2},$ $\mathfrak{a}_{0},$ \ $\mathfrak{a}_{-2}.$

\begin{eqnarray*}
\partial _{5}^{ij} &=&c_{k1}\mathfrak{a}_{2}+c_{k2}\mathfrak{a}_{0}+c_{k3}%
\mathfrak{a}_{-2}, \\
k &=&3,\text{ }i=7,9,10;\text{ }j=1,2,3 \\
&=&2,\text{ }i=13,15,16;\text{ }j=1,2,3 \\
&=&1,\text{ }i=23,25,26;\text{ }j=1,2,3 \\
\partial _{5}^{ij} &=&\left( 
\begin{array}{cc}
c_{\sigma 3} & c_{\tau 3} \\ 
c_{\sigma 2} & c_{\tau 2} \\ 
c_{\sigma 1} & c_{\tau 1}%
\end{array}%
\right) \\
\left( \sigma ,\tau \right) &=&\left( 2,3\right) ,\text{ }i=18,19,20;\text{ }%
j=4,5;\text{ signs }\underset{++}{\underset{--}{+-}} \\
&=&\left( 1,3\right) ,\text{ }i=28,29,30;\text{ }j=4,6;\text{ signs }%
\underset{-+}{\underset{+-}{-+}} \\
&=&\left( 1,2\right) ,\text{ }i=32,33,34;\text{ }j=5,6;\text{ signs }%
\underset{--}{\underset{++}{--}} \\
\partial _{5}^{ij} &=&0\text{ in other cases of }i,j.
\end{eqnarray*}

\begin{eqnarray*}
\partial _{4}^{ij} &=&\left( c_{\sigma 1}+c_{\tau 1}\right) \mathfrak{a}%
_{2}+\left( c_{\sigma 2}+c_{\tau 2}\right) \mathfrak{a}_{0}+\left( c_{\sigma
3}+c_{\tau 3}\right) \mathfrak{a}_{-2} \\
\left( \sigma ,\tau \right) &=&\left( 2,3\right) ,\text{ }i=12,14,15;\text{ }%
j=4,5,6;\text{ }j=7,8,9 \\
&=&\left( 1,3\right) ,\text{ }i=22,24,25;\text{ }j=4,5,6;\text{ }j=11,12,13
\\
&=&\left( 1,2\right) ,\text{ }i=28,30,31;\text{ }j=7,8,9;\text{ }j=11,12,13
\\
\partial _{4}^{ij} &=&\left( c_{kl}\right) _{k,l=1,2,3.}\text{ ; }i=5,9,19;%
\text{ }j=1,2,3\text{; signs\ }\underset{-\text{ \ }-\text{ \ }-}{\underset{+%
\text{ \ }+\text{ \ }+}{---}} \\
&=&\left( c_{kl}\right) _{k,l=1,2,3.}^{T}\text{ ; }i=33,34,35;\text{ }%
j=10,14,15\text{; signs\ }\underset{-\text{ \ }-\text{ \ }-}{\underset{+%
\text{ \ }+\text{ \ }+}{---}} \\
\partial _{4}^{ij} &=&0\text{ in another case of }i,j.
\end{eqnarray*}

\begin{eqnarray*}
\partial _{3}^{ij} &=&\left( c_{\sigma 1}+c_{\tau 1}+c_{\nu 1}\right) 
\mathfrak{a}_{2}+\left( c_{\sigma 2}+c_{\tau 2}+c_{\nu 2}\right) \mathfrak{a}%
_{0}+\left( c_{\sigma 3}+c_{\tau 3}+c_{\nu 3}\right) \mathfrak{a}_{-2} \\
\left( \sigma ,\tau ,\nu \right) &=&\left( 1,2,3\right) ,\text{ }%
i=18,20,21;j=7,8,9;j=13,14,15;j=16,17,18 \\
\partial _{3}^{ij} &=&\left( 
\begin{array}{ccc}
c_{\tau 1} & c_{\tau 2} & c_{\tau 3} \\ 
-c_{\sigma 1} & -c_{\sigma 2} & -c_{\sigma 3}%
\end{array}%
\right) \\
\left( \sigma ,\tau \right) &=&\left( 2,3\right) ,\text{ }i=11,15;j=10,11,12
\\
&=&\left( 1,3\right) ,\text{ }i=6,15;j=4,5,6 \\
&=&\left( 1,2\right) ,\text{ }i=6,11;j=1,2,3 \\
\partial _{3}^{ij} &=&0\text{ in other cases of }i,j
\end{eqnarray*}

\noindent
and finally

\begin{eqnarray*}
\partial _{2}^{ij} &=&\left( -c_{\sigma 1},-c_{\sigma 2},-c_{\sigma 3}\right)
\\
\sigma &=&1,\text{ }i=7,\text{ }j=1,2,3 \\
\sigma &=&2,\text{ }i=7,\text{ }j=4,5,6 \\
\sigma &=&3,\text{ }i=7,\text{ }j=7,8,9.
\end{eqnarray*}

The homology groups are in general $H_{i}\left( \mathfrak{a}_{E}\right)
=R^{\dim \ker \partial _{i}-\dim im\partial _{i+1}}.$ In this case

\begin{eqnarray*}
H_{7}\left( \mathfrak{a}_{E}\right) &=&R;\text{ }H_{6}\left( \mathfrak{a}%
_{E}\right) =R^{6} \\
H_{5}\left( \mathfrak{a}_{E}\right) &=&R^{20-rank\partial _{5}},\text{ }%
rank\partial _{5}=0,1,2,...,6 \\
H_{4}\left( \mathfrak{a}_{E}\right) &=&R^{35-rank\partial _{4}-rank\partial
_{5}},\text{ }rank\partial _{4}=0,1,2,...,15;rank\partial _{5}=0,1,2,...,6 \\
H_{3}\left( \mathfrak{a}_{E}\right) &=&R^{35-rank\partial _{3}-rank\partial
_{4}},\text{ }rank\partial _{3}=0,1,2,...,21;rank\partial _{4}=0,1,2,...,15
\\
H_{2}\left( \mathfrak{a}_{E}\right) &=&R^{20-rank\partial _{3}},\text{ }%
rank\partial _{3}=0,1,2,...,21 \\
H_{1}\left( \mathfrak{a}_{E}\right) &=&R^{6};\text{ }H_{0}\left( \mathfrak{a}%
_{E}\right) =R.
\end{eqnarray*}

The Euler characteristic is in every case $\mathcal{X}\left( \mathfrak{a}%
_{E}\right) =0.$

$\mathit{\square }$

{}


\medskip

\begin{thebibliography}{9}

\bibitem{Baez-Morton} 
Baez,J., Morton and Vicary 
{\it On the Categorified Heisenberg Algebra.}
http://golem.ph.utexas.edu/category/2012/07/morton\_and\_vicary\_on\_the%
\_categ.html

\bibitem{Blasco} 
Blasco, X., Mart\'{\i}nez, M.,  Herrero, J.M., Ramos, C.,
Sanchis, J.  
{\it Model-based predictive control of greenhouse climate for
reducing energy  and water consumption.} 
Computers and 
 Electronics in
Agriculture 55. (2007) 49-70.

\bibitem{Cairns} 
 Cairns, Grant; Jambor, Sebastian  
{\it The cohomology of the Heinsenberg Lie algebras over fields of finite characteristic.}
Proc. Amer. Math. Soc. 136 (2008), no. 11, 3803--3807. 

\bibitem{Merlin Oliver} 
Merlin Olivier, Ahmad Al Bitar, Vincent Rivalland, Pierre B\'{e}%
ziat, Eric Ceschia, G\'{e}rard Dedieu,  
{\it \textquotedblleft An Analytical Model of Evaporation Efficiency for 
Unsaturated Soil Surfaces
with an Arbitrary Thickness.}
\textquotedblright\ J. Appl. Meteor. Climatol. 50, (2011)
 457--471.

\bibitem{Morton} 
Morton, J. C. and Vicary, J. 
{\it The Categorified Heisenberg Algebra I:
A Combinatorial Representation.} arXiv:1207.2054v2 [math.QA].

\bibitem{Sanchez} 
 S\'anchez-Gabites, J. J. 
 {\it Dynamical systems and shapes.} Rev. R. Acad. Cienc. Exactas F\'is. Nat. Ser. A Math. RACSAM 102 (2008), no. 1, 
 127--159.

\bibitem{Santaroubane} 
Santaroubane, L. J.  
{\it Cohomology of Heinsenberg Lie algebras.}
 Proceedings of the American Mathematical Society 87 (1983), no. 1, 23--28.

\bibitem{Tkachuk} 
Tkachuk, V. M. 
{\it Deformed Heisenberg algebra with minimal length
and equivalence principle.} arXiv:1301.1891v1 [gr-qc].
\end{thebibliography}
\end{document}